\documentclass[12pt,vatola]{article}

\usepackage{graphicx}

\def\C{\centerline}

\def\re#1{\par\hangindent\parindent\indent\llap{#1\enspace}\ignorespaces}

\def\no{\noindent}

\begin{document}

\vskip 15mm

\C{\large\bf On Algebraic Multi-Ring Spaces  }  \vskip 5mm

\C{Linfan Mao} \vskip 3mm \C{\scriptsize (Academy of Mathematics
and System Sciences, Chinese Academy of Sciences, Beijing 100080)}

\vskip 8mm
\begin{minipage}{130mm}
\no{\bf Abstract}: {\small A Smarandache multi-space is a union of
$n$ spaces $A_1,A_2,\cdots ,A_n$ with some additional conditions
holding. Combining Smarandache multi-spaces with rings in
classical ring theory, the conception of multi-ring spaces is
introduced. Some characteristics of a multi-ring space are
obtained in this paper}

\vskip 2mm \no{\bf Key words:} {\small  ring, multi-space,
multi-ring space, ideal subspace chain.}

 \vskip 2mm \no{{\bf
Classification:} AMS(2000) 16U80, 16W25}
\end{minipage}

\vskip 8mm

{\bf $1$. Introduction}

\vskip 6mm

The notion of multi-spaces is introduced by Smarandache in $[6]$
under his idea of hybrid mathematics: {\it  combining different
fields into a unifying field}($[7]$), which is defined as follows.

\vskip 4mm

\no{\bf Definition $1.1$} \ {\it For any integer $i, 1\leq i\leq
n$ let $A_i$ be a set with ensemble of law $L_i$, and the
intersection of $k$ sets $A_{i_1},A_{i_2},\cdots , A_{i_k}$ of
them constrains the law $I(A_{i_1},A_{i_2},\cdots , A_{i_k})$.
Then the union of $A_i$, $1\leq i\leq n$

$$\widetilde{A} \ = \ \bigcup\limits_{i=1}^n A_i$$

\no is called a multi-space.}

\vskip 3mm

As we known, a set $R$ with two binary operation ¡°$+$¡± and
¡°$\circ$¡±, denoted by $(R \ ;+,\circ )$, is said to be a {\it
ring} if for $\forall x,y\in R$, $x+y\in R$, $x\circ y\in R$, the
following conditions hold.

\vskip 3mm

($i$) \ $(R \ ;+)$ is an abelian group;

($ii$) \ $(R \ ;\circ )$ is a semigroup;

($iii$) \ For $\forall x,y,z\in R$, $x\circ (y+z)=x\circ y+x\circ
z$ and $(x+y)\circ z= x\circ z+y\circ z.$

By combining Smarandache multi-spaces with rings, a new kind of
algebraic structure called multi-ring space is found, which is
defined in the following.

\vskip 4mm

\no{\bf Definition $1.2$} \ {\it Let
$\widetilde{R}=\bigcup\limits_{i=1}^mR_i$ be a complete
multi-space with double binary operation set
$O(\widetilde{R})=\{(+_i,\times_i) , 1\leq i\leq m\}$. If for any
integers $i, j, \ i\not= j, 1\leq i, j\leq m$, $(R_i; +_i,
\times_i)$ is a ring and for $\forall x,y,z\in\widetilde{R}$,

$$ (x+_iy)+_jz = x+_i(y+_jz), \  \ \ (x\times_iy)\times_jz = x\times_i(y\times_jz)$$

\no and

$$x\times_i(y+_jz) = x\times_iy +_jx\times_iz, \  \ \ (y+_jz)\times_ix = y\times_ix +_jz\times_ix$$

\no if all their operation results exist, then $\widetilde{R}$ is
called a multi-ring space. If for any integer $1\leq i\leq m$,
$(R;+_i,\times_i)$ is a filed, then $\widetilde{R}$ is called a
multi-filed space.}

\vskip 3mm

For a multi-ring space $\widetilde{R}=\bigcup\limits_{i=1}^mR_i$,
let $\widetilde{S}\subset\widetilde{R}$ and
$O(\widetilde{S})\subset O(\widetilde{R})$, if $\widetilde{S}$ is
also a multi-ring space with double binary operation set
$O(\widetilde{S})$ , then call $\widetilde{S}$ a {\it multi-ring
subspace} of $\widetilde{R}$. We have the following criterions for
the multi-ring subspaces.

The subject of this paper is to find some characteristics of a
multi-ring space. For terminology and notation not defined here
can be seen in $[1],[5],[12]$ for algebraic terminologies and in
$[2],[6]-[11]$ for multi-spaces and logics.

\vskip 8mm

{\bf $2.$ Characteristics of a multi-ring space}

\vskip 5mm

First, we have the following result for multi-ring subspace of a
multi-ring space.

\vskip 4mm

\no{\bf Theorem $2.1$} \ {\it For a multi-ring space
$\widetilde{R}=\bigcup\limits_{i=1}^mR_i$, a subset
$\widetilde{S}\subset\widetilde{R}$ with $O(\widetilde{S})\subset
O(\widetilde{R})$ is a multi-ring subspace of $\widetilde{R}$ if
and only if for any integer $k, 1\leq k\leq m$,
$(\widetilde{S}\bigcap R_k; +_k,\times_k)$ is a subring of $(R_k;
+_k, \times_k)$ or $\widetilde{S}\bigcap R_k=\emptyset$.}

\vskip 3mm

{\it Proof} \ For any integer $k, 1\leq k\leq m$, if
$(\widetilde{S}\bigcap R_k; +_k,\times_k)$ is a subring of $(R_k;
+_k,\times_k)$ or $\widetilde{S}\bigcap R_k=\emptyset$, then since
$\widetilde{S}=\bigcup\limits_{i=1}^m(\widetilde{S}\bigcap R_i)$,
we know that $\widetilde{S}$ is a multi-ring subspace by
definition of a multi-ring space.

Now if $\widetilde{S}=\bigcup\limits_{j=1}^sS_{i_j}$ is a
multi-ring subspace of $\widetilde{R}$ with double binary
operation set $O(\widetilde{S})=\{(+_{i_j},\times_{i_j}), 1\leq
j\leq s\}$, then $(S_{i_j};  +_{i_j},\times_{i_j})$ is a subring
of $(R_{i_j}; +_{i_j}, \times_{i_j})$. Therefore, for any integer
$j, 1\leq j\leq s$, $S_{i_j}=R_{i_j}\bigcap\widetilde{S}$. But for
other integer $l\in\{i; 1\leq i\leq m\}\setminus\{i_j; 1\leq j\leq
s\}$, $\widetilde{S}\bigcap S_l=\emptyset$. \quad\quad $\natural$

Applying the criterions for subrings of a ring, we get the
following result.

\vskip 4mm

\no{\bf Theorem $2.2$} \ {\it For a multi-ring space
$\widetilde{R}=\bigcup\limits_{i=1}^mR_i$, a subset
$\widetilde{S}\subset\widetilde{R}$ with $O(\widetilde{S})\subset
O(\widetilde{R})$ is a multi-ring subspace of $\widetilde{R}$ if
and only if for any double binary operations $(+_j, \times_j)\in
O(\widetilde{S})$, $(\widetilde{S}\bigcap R_j; +_j) \prec
(R_j;+_j)$ and $(\widetilde{S}; \times_j)$ is complete.}

\vskip 3mm

{\it Proof} \ According to Theorem $2.1$, we know that
$\widetilde{S}$ is a multi-ring subspace if and only if for any
integer $i, 1\leq i\leq m$, $(\widetilde{S}\bigcap
R_i;+_i,\times_i)$ is a subring of $(R_i;+_i,\times_i)$ or
$\widetilde{S}\bigcap R_i=\emptyset$. By a well known criterions
for subrings of a ring (see also $[5]$), we know that
$(\widetilde{S}\bigcap R_i;+_i,\times_i)$ is a subring of
$(R_i;+_i,\times_i)$ if and only if for any double binary
operations $(+_j ,\times_j)\in O(\widetilde{S})$,
$(\widetilde{S}\bigcap R_j; +_j) \prec (R_j;+_j)$ and
$(\widetilde{S}; \times_j)$ is a complete set. This completes the
proof. \quad\quad $\natural$

We use the {\it ideal subspace chain} of a multi-ring space to
characteristic its structure properties. An {\it ideal subspace}
$\widetilde{I}$ of a multi-ring space
$\widetilde{R}=\bigcup\limits_{i=1}^mR_i$ with double binary
operation set $O(\widetilde{R})$ is a multi-ring subspace of
$\widetilde{R}$ satisfying the following conditions:

\vskip 3mm

$(i)$ \ $\widetilde{I}$ is a multi-group subspace with operation
set $\{+ | \ (+ ,\times)\in O(\widetilde{I})\}$;

$(ii)$ \ for any $r\in\widetilde{R}, a\in\widetilde{I}$ and
$(+,\times)\in O(\widetilde{I})$, $r\times a\in\widetilde{I}$ and
$a\times r\in\widetilde{I}$ if their operation results exist.

\vskip 4mm

\no{\bf Theorem $2.3$} \ {\it A subset $\widetilde{I}$ with
$O(\widetilde{I}), O(\widetilde{I})\subset O(\widetilde{R})$ of a
multi-ring space $\widetilde{R}=\bigcup\limits_{i=1}^mR_i$ with
double binary operation set $O(\widetilde{R})=\{(+_i,\times_i)| \
1\leq i\leq m\}$ is an ideal subspace if and only if for any
integer $i, 1\leq i\leq m$, $(\widetilde{I}\bigcap R_i,
+_i,\times_i)$ is an ideal of the ring $(R_i,+_i,\times_i)$ or
$\widetilde{I}\bigcap R_i=\emptyset$.}

\vskip 3mm

{\it Proof} \ By definition of an ideal subspace, the necessity of
the condition is obvious.

For the sufficiency, denote by $\widetilde{R}(+,\times)$ the set
of elements in $\widetilde{R}$ with binary operations ¡°$+$¡± and
¡°$\times$¡±. If there exists an integer $i$ such that
$\widetilde{I}\bigcap R_i\not=\emptyset$ and
$(\widetilde{I}\bigcap R_i, +_i,\times_i)$ is an ideal of
$(R_i,+_i,\times_i)$, then for $\forall a\in\widetilde{I}\bigcap
R_i$, $\forall r_i\in R_i$, we know that

$$r_i\times_i a\in\widetilde{I}\bigcap R_i; \ \ \
a\times_i r_i\in\widetilde{I}\bigcap R_i .$$

Notice that $\widetilde{R}(+_i,\times_i)=R_i$. Therefore, we get
that for $\forall r\in\widetilde{R}$,

$$r\times_i a\in\widetilde{I}\bigcap R_i; \ {\rm and} \
a\times_i r\in\widetilde{I}\bigcap R_i ,$$

\no if their operation result exist. Whence, $\widetilde{I}$ is an
ideal subspace of $\widetilde{R}$. \quad\quad $\natural$

An ideal subspace $\widetilde{I}$ of a multi-ring space
$\widetilde{R}$ is said {\it maximal} if for any ideal subspace
$\widetilde{I}'$, if
$\widetilde{R}\supseteq\widetilde{I}'\supseteq\widetilde{I}$, then
$\widetilde{I}'=\widetilde{R}$ or $\widetilde{I}'=\widetilde{I}$.
For any order of the double binary operation set
$O(\widetilde{R})$ of a multi-ring space
$\widetilde{R}=\bigcup\limits_{i=1}^mR_i$, not loss of generality,
assume it being $(+_1,\times_1)\succ (+_2,\times_2)\succ
\cdots\succ(+_m,\times_m)$, we can define an {\it ideal subspace
chain} of $\widetilde{R}$ by the following programming.

\vskip 3mm

$(i)$ Construct the ideal subspace chain

$$\widetilde{R}\supset\widetilde{R}_{11}\supset\widetilde{R}_{12}
\supset\cdots\supset\widetilde{R}_{1s_1}$$

\no under the double binary operation $(+_1,\times_1)$, where
$\widetilde{R}_{11}$ is a maximal ideal subspace of
$\widetilde{R}$ and in general, for any integer $i, \ 1\leq i\leq
m-1$, $\widetilde{R}_{1(i+1)}$ is a maximal ideal subspace of
$\widetilde{R}_{1i}$.

$(ii)$ If the ideal subspace

$$\widetilde{R}\supset\widetilde{R}_{11}\supset\widetilde{R}_{12}
\supset\cdots\supset\widetilde{R}_{1s_1}\supset\cdots\supset\widetilde{R}_{i1}
\supset\cdots\supset\widetilde{R}_{is_i}$$

\no has been constructed for $(+_1,\times_1)\succ
(+_2,\times_2)\succ \cdots\succ(+_i,\times_i)$, $1\leq i\leq m-1$,
then construct an ideal subspace chain of $\widetilde{R}_{is_i}$

$$\widetilde{R}_{is_i}\supset\widetilde{R}_{(i+1)1}\supset\widetilde{R}_{(i+1)2}
\supset\cdots\supset\widetilde{R}_{(i+1)s_1}$$

\no under the operations $(+_{i+1},\times_{i+1})$,  where
$\widetilde{R}_{(i+1)1}$ is a maximal ideal subspace of
$\widetilde{R}_{is_i}$ and in general,
$\widetilde{R}_{(i+1)(i+1)}$ is a maximal ideal subspace of
$\widetilde{R}_{(i+1)j}$ for any integer $j, 1\leq j\leq s_i-1$.
Define the ideal subspace chain of $\widetilde{R}$ under
$(+_1,\times_1)\succ (+_2,\times_2)\succ
\cdots\succ(+_{i+1},\times_{i+1})$ being

$$\widetilde{R}\supset\widetilde{R}_{11}
\supset\cdots\supset\widetilde{R}_{1s_1}\supset\cdots\supset\widetilde{R}_{i1}
\supset\cdots\supset\widetilde{R}_{is_i}\supset\widetilde{R}_{(i+1)1}
\supset\cdots\supset\widetilde{R}_{(i+1)s_{i+1}}.$$

Similar to a multi-group space([$3$]), we have the following
results for the ideal subspace chain of a multi-ring space.

\vskip 4mm

\no{\bf Theorem $2.4$} \ {\it For a multi-ring space
$\widetilde{R}=\bigcup\limits_{i=1}^mR_i$, its ideal subspace
chain only has finite terms if and only if for any integer $i,
1\leq i\leq m$, the ideal chain of the ring $(R_i;+_i,\times_i)$
has finite terms, i.e., each ring $(R_i;+_i,\times_i)$ is an Artin
ring.}

\vskip 3mm

{\it Proof} \ Let the order of double operations in
$\overrightarrow{O}(\widetilde{R})$ be

$$(+_1,\times_1)\succ (+_2,\times_2)\succ\cdots\succ (+_m,\times_m)$$

\no and a maximal ideal chain in the ring $(R_1;+_1,\times_1)$ is

$$R_1\succ R_{11}\succ\cdots\succ R_{1t_1}.$$

\no Calculate

$$\widetilde{R}_{11}=\widetilde{R}\setminus\{R_1\setminus R_{11}\}=
R_{11}\bigcup (\bigcup\limits_{i=2}^m)R_i,$$

$$\widetilde{R}_{12}=\widetilde{R}_{11}\setminus\{R_{11}\setminus R_{12}\}=
R_{12}\bigcup (\bigcup\limits_{i=2}^m)R_i,$$

$$\cdots\cdots\cdots\cdots\cdots\cdots$$

$$\widetilde{R}_{1t_1}=\widetilde{R}_{1t_1}\setminus\{R_{1(t_1-1)}\setminus R_{1t_1}\}=
R_{1t_1}\bigcup (\bigcup\limits_{i=2}^m)R_i.$$

\no According to Theorem $3.10$, we know that

$$\widetilde{R}\supset\widetilde{R}_{11}\supset\widetilde{R}_{12}
\supset\cdots\supset\widetilde{R}_{1t_1}$$

\no  is a maximal ideal subspace chain of $\widetilde{R}$ under
the double binary operation $(+_1,\times_1)$. In general, for any
integer $i, 1\leq i\leq m-1$, assume

$$R_i\succ R_{i1}\succ\cdots\succ R_{it_i}$$

\no is a maximal ideal chain in the ring
$(R_{(i-1)t_{i-1}};+_i,\times_i)$. Calculate

$$\widetilde{R}_{ik}=R_{ik}\bigcup(\bigcup\limits_{j=i+1}^m)\widetilde{R}_{ik}\bigcap R_i$$

\no Then we know that

$$\widetilde{R}_{(i-1)t_{i-1}}\supset\widetilde{R}_{i1}\supset\widetilde{R}_{i2}
\supset\cdots\supset\widetilde{R}_{it_i}$$

\no is a maximal ideal subspace chain of
$\widetilde{R}_{(i-1)t_{i-1}}$ under the double operation
$(+_i,\times_i)$ by Theorem $3.10$. Whence, if for any integer $i,
1\leq i\leq m$, the ideal chain of the ring $(R_i;+_i,\times_i)$
has finite terms, then the ideal subspace chain of the multi-ring
space $\widetilde{R}$ only has finite terms and if there exists
one integer $i_0$ such that the ideal chain of the ring
$(R_{i_0},+_{i_0},\times_{i_0})$ has infinite terms, then there
must be infinite terms in the ideal subspace chain of the
multi-ring space $\widetilde{R}$. \quad\quad $\natural$.

A multi-ring space is called an {\it Artin multi-ring space} if
each ideal subspace chain only has finite terms. We have the
following corollary by Theorem $3.11$.

\vskip 4mm

\no{\bf Corollary $2.1$} {\it A multi-ring space
$\widetilde{R}=\bigcup\limits_{i=1}^m$ with double binary
operation set $O(\widetilde{R})=\{(+_i,\times_i)| \ 1\leq i\leq
m\}$ is an Artin multi-ring space if and only if for any integer
$i, 1\leq i\leq m$, the ring $(R_i;+_i,\times_i)$ is an Artin
ring.}

\vskip 3mm

For a multi-ring space $\widetilde{R}=\bigcup\limits_{i=1}^m$ with
double binary operation set $O(\widetilde{R})=\{(+_i,\times_i)| \
1\leq i\leq m\}$, an element $e$ is an {\it idempotent} element if
$e_{\times}^2 = e\times e = e$ for a double binary operation
$(+,\times)\in O(\widetilde{R})$. We define the {\it directed sum}
$\widetilde{I}$ of two ideal subspaces $\widetilde{I}_1$ and
$\widetilde{I}_2$ as follows:

\vskip 3mm

$(i)$ $\widetilde{I}=\widetilde{I}_1\bigcup\widetilde{I}_2$;

$(ii)$ $\widetilde{I}_1\bigcap\widetilde{I}_2=\{0_+\}, \ {\rm or}
\ \widetilde{I}_1\bigcap\widetilde{I}_2=\emptyset$, where $0_+$
denotes an unit element under the operation $+$.

Denote the directed sum of $\widetilde{I}_1$ and $\widetilde{I}_2$
by

$$\widetilde{I}=\widetilde{I}_1\bigoplus\widetilde{I}_2.$$

If for any $\widetilde{I}_1, \widetilde{I}_2$,
$\widetilde{I}=\widetilde{I}_1\bigoplus\widetilde{I}_2$ implies
that $\widetilde{I}_1=\widetilde{I}$ or
$\widetilde{I}_2=\widetilde{I}$, then $\widetilde{I}$ is called
{\it non-reducible}. We have the following result for the Artin
multi-ring space similar to a well-known result for the Artin ring
(see $[12]$).

\vskip 4mm

\no{\bf Theorem $2.5$} \ {\it Any Artin multi-ring space
$\widetilde{R}=\bigcup\limits_{i=1}^mR_i$ with double binary
operation set $O(\widetilde{R})=\{(+_i,\times_i)| \ 1\leq i\leq
m\}$ is a directed sum of finite non-reducible ideal subspaces,
and if for any integer $i,1\leq i\leq m$, $(R_i;+_i,\times_i)$ has
unit $1_{\times_i}$, then}

$$\widetilde{R}=\bigoplus\limits_{i=1}^m(\bigoplus\limits_{j=1}^{s_i}(R_i\times_ie_{ij})
\bigcup (e_{ij}\times_i R_i)),$$

\no{\it where $e_{ij}, 1\leq j\leq s_i$ are orthogonal idempotent
elements of the ring $R_i$.}

\vskip 3mm

{\it Proof} \ Denote by $\widetilde{M}$ the set of ideal subspaces
which can not be represented by a directed sum of finite ideal
subspaces in $\widetilde{R}$. According to Theorem $3.11$, there
is a minimal ideal subspace $\widetilde{I}_0$ in $\widetilde{M}$.
It is obvious that $\widetilde{I}_0$ is reducible.

Assume that $\widetilde{I}_0=\widetilde{I}_1+\widetilde{I}_2$.
Then $\widetilde{I}_1\not\in\widetilde{M}$ and
$\widetilde{I}_2\not\in\widetilde{M}$. Therefore,
$\widetilde{I}_1$ and $\widetilde{I}_2$ can be represented by
directed sums of finite ideal subspaces. Whence, $\widetilde{I}_0$
can be also represented  by a directed sum of finite ideal
subspaces. Contradicts that $\widetilde{I}_0\in\widetilde{M}$.

Now let

$$\widetilde{R}=\bigoplus\limits_{i=1}^s\widetilde{I}_i,$$

\no where each $\widetilde{I}_i, 1\leq i\leq s$, is non-reducible.
Notice that for a double operation $(+,\times )$, each
non-reducible ideal subspace of $\widetilde{R}$ has the form

$$(e\times R(\times))\bigcup (R(\times )\times e), \ \ e\in R(\times ).$$

\no Whence, we know that there is a set $T\subset\widetilde{R}$
such that

$$\widetilde{R}=\bigoplus\limits_{e\in T, \ \times\in O(\widetilde{R})}
(e\times R(\times))\bigcup (R(\times )\times e).$$

For any operation $\times\in O(\widetilde{R})$ and the unit
$1_{\times}$, assume that

$$1_{\times} = e_1\oplus e_2\oplus\cdots \oplus e_l, \ e_i\in T, \ 1\leq i\leq s.$$

\no Then

$$e_i\times 1_{\times}= (e_i\times e_1)\oplus (e_i\times e_2)\oplus\cdots \oplus (e_i\times e_l).$$

\no Therefore, we get that

$$e_i = e_i\times e_i=e_i^2 \ \ {\rm and} \ \ e_i\times e_j=0_i \ \ {\rm for} \ \ i\not=j.$$

\no That is, $e_i, 1\leq i\leq l$, are orthogonal idempotent
elements of $\widetilde{R}(\times)$. Notice that
$\widetilde{R}(\times)=R_h$ for some integer $h$. We know that
$e_i, 1\leq i\leq l$ are orthogonal idempotent elements of the
ring $(R_h,+_h,\times_h)$. Denoted by $e_{hj}$ for $e_j$, $1\leq
j\leq l$. Consider all units in $\widetilde{R}$, we get that

$$\widetilde{R}=\bigoplus\limits_{i=1}^m(\bigoplus\limits_{j=1}^{s_i}(R_i\times_ie_{ij})
\bigcup (e_{ij}\times_i R_i)).$$

\no This completes the proof. \quad\quad $\natural$

\vskip 4mm

\no{\bf Corollary $2.2$}([12]) \ {\it Any Artin ring $(R \
;+,\times )$ is a directed sum of finite ideals, and if $(R \
;+,\times )$ has unit $1_{\times}$, then}

$$R =\bigoplus\limits_{i=1}^s R_ie_{i },$$

\no{\it where $e_{i}, 1\leq i\leq s$ are orthogonal idempotent
elements of the ring $(R;+,\times )$.}

\vskip 8mm

{\bf $3.$ Open problems for a multi-ring space}

\vskip 5mm

Similar to Artin multi-ring space, we can also define {\it Noether
multi-ring spaces, simple multi-ring spaces, half-simple
multi-ring spaces, $\cdots$, etc.}.  The open problems for these
new algebraic structure are as follows.

\vskip 3mm

\no{\bf Problem $3.1$} \ Call a ring $R$ a Noether ring if its
every ideal chain only has finite terms. Similarly, for a
multi-ring space $\widetilde{R}$, if its every ideal multi-ring
subspace chain only has finite terms, it is called a Noether
multi-ring space.  {\it Whether can we find its structures similar
to Corollary $2.2$ and Theorem $2.5$?}

\vskip 3mm

\no{\bf Problem $3.2$} \ {\it Similar to ring theory, define a
Jacobson or Brown-McCoy radical for multi-ring spaces and
determine their contribution to multi-ring spaces.}

\vskip 8mm

{\bf References}

\vskip 5mm

\re{[1]}G.Birkhoff and S.Mac Lane, {\it A Survey of Modern
Algebra}, Macmillan Publishing Co., Inc, 1977.

\re{[2]}Daniel Deleanu, {\it A Dictionary of Smarandache
Mathematics}, Buxton University Press, London \& New York,2004.

\re{[3]}L.F.Mao, On Algebraic Multi-Group Spaces, {\it eprint
arXiv: math/0510427}, 10/2005.

\re{[4]}L.F.Mao, {\it Automorphism Groups of Maps, Surfaces and
Smarandache Geometries}, American Research Press, 2005.

\re{[5]}L.Z. Nie and S.S Ding, {\it Introduction to Algebra},
Higher Education Publishing Press, 1994.

\re{[6]} F.Smarandache, Mixed noneuclidean geometries, {\it eprint
arXiv: math/0010119}, 10/2000.

\re{[7]}F.Smarandache, {\it A Unifying Field in Logics.
Neutrosopy: Neturosophic Probability, Set, and Logic}, American
research Press, Rehoboth, 1999.

\re{[8]}F.Smarandache, Neutrosophy, a new Branch of Philosophy,
{\it Multi-Valued Logic}, Vol.8, No.3(2002)(special issue on
Neutrosophy and Neutrosophic Logic), 297-384.

\re{[9]}F.Smarandache, A Unifying Field in Logic: Neutrosophic
Field, {\it Multi-Valued Logic}, Vol.8, No.3(2002)(special issue
on Neutrosophy and Neutrosophic Logic), 385-438.

\re{[10]]}W.B.Vasantha Kandasamy, {\it Bialgebraic structures and
Smarandache bialgebraic structures}, American Research Press,
2003.

\re{[11]}W.B.Vasantha Kandasamy and F.Smarandache, {\it Basic
Neutrosophic Algebraic Structures and Their Applications to Fuzzy
and Neutrosophic Models}, HEXIS, Church Rock, 2004.

\re{[12]}Quanyan Xong, {\it Ring Theory}, Wuhan University Press,
1993.

\end{document}